\documentclass[12pt]{amsart}
\usepackage{amsfonts,latexsym,amsthm,amssymb}
\usepackage[all]{xy}

\input epsf
\CompileMatrices

\newtheorem{theorem}{Theorem}[section]
\newtheorem{lemma}[theorem]{Lemma}
\newtheorem{corol}[theorem]{Corollary}
\newtheorem{prop}[theorem]{Proposition}

{\theoremstyle{definition} \newtheorem{defin}[theorem]{Definition}}
{\theoremstyle{remark} \newtheorem{remark}[theorem]{Remark}
\newtheorem{example}[theorem]{Example}}

\newcommand{\Pbb}{{\mathbb{P}}}

\newcommand{\Til}[1]{{\widetilde{#1}}}

\newcommand{\csm}{c_{\text{\rm SM}}}
\newcommand{\ssm}{s^\circ}

\title{Inclusion-exclusion and Segre classes, II}
\author{Paolo Aluffi}
\address{Max-Planck-Institut f\"ur Mathematik, Bonn, Germany}
\address{Florida State University, Tallahassee, Florida}

\makeindex

\begin{document}

\begin{abstract}
Considerations based on the known relation between different
characteristic classes for singular hypersufaces suggest that a form
of the `inclusion-exclusion' principle may hold for Segre classes. We
formulate and prove such a principle for a notion closely related to
Segre classes. This is used to provide a simple computation of the
classes introduced in \cite{incexcI}, in certain special (but
representative) cases.
\end{abstract}

\maketitle


\section{Introduction}

Recent work on relations between the Chern-Schwartz-MacPherson class
of a hypersurface $X$ in a nonsingular ambient variety $M$ and the
class of its virtual tangent bundle has revealed a tight connection
between the former and the {\em Segre class\/} of the Jacobian
subscheme of $X$ (cf.~\cite{MR96d:14004} and \cite{MR2001i:14009},
\S1.1). A strong motivation to pursue this connection is the important
r\^ole played by Segre classes in intersection theory---the hope is
that functoriality properties of Segre classes (such as those that may
be inherited via the connection with Chern-Schwartz-MacPherson's
classes) would lead to new computational tools for Segre classes, or
at least point in the right direction to look for such tools. In
this paper we discuss an `inclusion-exclusion principle' for Segre
classes, inspired by this connection.

This article is a counterpoint to \cite{incexcI}, where we have {\em
  imposed\/} an inclusion-exclusion principle on a Segre class-type
  notion. The resulting {\em SM-Segre class\/} satisfies a number of
  remarkable properties, which must be a reflection of unknown and
  potentially useful properties of ordinary Segre classes.
This leads us in this article to search for other instances where
an inclusion-exclusion principle may be at work in the theory of Segre
classes.

By {\em inclusion-exclusion\/} we refer to the familiar counting
principle according to which the number of elements in the
intersection of a family of finite sets may be computed by adding the
cardinalities of the sets, subtracting the cardinalities of their
pairwise unions, adding back the cardinalities of triple unions,
etc. An analog of this principle is trivially satisfied by the
topological Euler characteristic, and the simplest form of the
functoriality property of the Chern-Schwartz-MacPherson is an
expression of the same principle. This is the observation leading to
the definition of SM-Segre classes in \cite{incexcI}.

We can now abstract away from Chern-Schwartz-MacPherson classes
for a moment, and look for other situations where Segre classes
express a behavior reminiscent of inclusion-exclusion. Our candidate
in this paper is proposed in \S\ref{candid}. We show that with a
suitable definition of the {\em union Segre class,\/} which
unfortunately is {\em not\/} quite the Segre class of the union, a
straightforward inclusion-exclusion formula computes the
(conventional) Segre class of the intersection $Y$ of closed
subschemes $X_1,\dots,X_n$ in an ambient irreducible scheme~$M$. In
fact we give a rather broad statement (Theorem~\ref{main}), which
specializes to inclusion-exclusion but encompasses a substantially
more general situation; indeed, some information can be obtained as
soon as $Y$ is contained in the intersection of the $X_i$, provided
(maybe surprisingly) that {\em enough\/} $X_i$ are considered.

After the fact, we go back to Chern-Schwartz-MacPherson classes: in
our main application of Theorem~\ref{main}, we observe that the union
Segre class in fact equals the SM-Segre class studied in
\cite{incexcI}, in a particularly well-behaved class of examples
(Theorem~\ref{relSM}, Corollary~\ref{surpr}). This recovers
immediately the main result of \cite{MR2001i:14009} in a very
particular, but representative case (cf.~Example~\ref{repr}). Our
hope, and our main motivation in this paper, is that applying
Theorem~\ref{relSM} to similarly representative situations may suggest
how to extend \cite{MR2001i:14009} to a wider class of varieties,
e.g., complete intersections. However, this issue is not further
explored here.

We also include in \S\ref{SM} comments meant to illustrate
the inclusion-exclusion principle in action as the number of loci
$X_i$ increases, and a few explicit examples.
\vskip 6pt
\noindent{\bf Acknowledgments} I thank the Max-Planck-Institut f\"ur
Mathematik in Bonn, Germany, for the hospitality and support.


\section{Union Segre classes and the inclusion-exclusion formula}
\label{candid}

We work in a fixed ambient variety $M$.

As we are aiming for an inclusion-exclusion formula, we first recall
the familiar counting analog for finite sets. Assume a set $Y$ is the
intersection of a finite number of finite sets $X_i$:
$$Y=X_1\cap \dots\cap X_r\quad.$$
Let $X_{i_1\dots i_s}$ denote the union $X_{i_1}\cup\dots\cup
X_{i_s}$; then the cardinality of $Y$ can be computed in terms of the
cardinalities of these unions by the formula
$$\# Y=\sum_{s=1}^r (-1)^{s-1} \sum_{i_1<\dots<i_s} \# X_{i_1\dots
i_s}\quad,$$
as is proved immediately by induction on $r$.

Using the same notations in the algebro--geometric situation, assume
$Y$ and the $X_i$ are proper subschemes of the ambient $M$; then we
can aim to a similar formula for the Segre class $s(Y,M)$:
\begin{equation}\label{segreform}
s(Y,M)=\sum_{s=1}^r (-1)^{s-1} \sum_{i_1<\dots<i_s} s(Y;X_{i_1},\dots,
X_{i_s};M)\quad;
\end{equation}
such an equality could be sought, for example, in the Chow group of
the union $\cup X_i$ (and we will be omitting evident push-forward
notations). The question we explore is how to define the class
$$s(Y;X_{i_1},\dots,X_{i_s};M)$$
so that such a formula may hold.

\begin{remark}
Counter to the first possible guess, this should not be the class of
the `union': for example, take $Y=p$, a point in $M=\Pbb^2$, and $X_1$,
$X_2$ two lines meeting at the point; then
$$s(X_1,\Pbb^2)+s(X_2,\Pbb^2)-s(X_1\cup X_2,\Pbb^2)=2[p]$$
rather than $[p]$.
\end{remark}

On the basis of a few such simple examples, one may reach the
conclusion that we have no right to expect a formula such as
(\ref{segreform}) to hold.  The surprise is that on the contrary there
{\em is\/} a straightforward way to define a Segre--like class
$s(Y;X_{i_1},\dots,X_{i_s};M)$ so that inclusion--exclusion works;
this class is defined on $X_{i_1}\cup\dots\cup X_{i_s}$, as it should,
but it is not (defined as) the Segre class of a subscheme in $M$
(unless $s=1$, cf.~Proposition~\ref{elem} below).

In fact, the class will depend on the specific choice of
$X_{i_1},\dots,X_{i_s}$ and on $Y$, and not only on the scheme
$X_{i_1}\cup\dots\cup X_{i_s}$; this seems an undesirable feature, but
we don't see any way around it at this moment. On the other hand, we
will see in \S\ref{SM} that the class is effectively computable
in several representative situations.

We will give the definition for proper closed subschemes
$X_1,\dots,X_s$ of $M$, all containing a closed subscheme $Y$ in fact
this requirement could be dropped, but at the expense of any geometric
interpretation of the classes arising in the process.
Let $\pi:\Til M \to M$ be a proper birational morphism such that $Y$
and all the $X_i$ pull--back to Cartier divisors
$\overline Y,\overline X_1,\dots,\overline X_s$ in $\Til M$; we say
then that $\pi$ is a `resolving' morphism. By construction, all the
$\overline X_i$ contain $\overline Y$ as a component, so that we can write 
$$\overline X_i=\overline Y+R_i$$ 
for well--defined residual divisors $R_i$ in $\Til M$.

\begin{defin}\label{usc}
With notations as above, we define the {\em union Segre class\/} of
$X_1,\dots,X_s$ in $M$, w.r.t.~$Y$, to be the class
$$s(Y;X_1,\dots,X_s;M)=\pi_*\frac{[R_1+\dots+R_s+\overline Y]}
{1+R_1+\dots+R_s+\overline Y}$$
in $A_*(X_1\cup\dots\cup X_s)$.
\end{defin}

Of course a resolving morphism $\pi$ exists (blow-up everything in
sight). However we must prove that the result is independent of the
choice of $\pi$.

\begin{lemma}
The definition of $s(Y;X_1,\dots,X_s;M)$ is independent of the chosen
resolving morphism $\pi$.
\end{lemma}

\begin{proof}
If $\pi:\Til M\to M$, $\pi':\Til M'\to M$ are two proper birational
morphisms for which $Y$ and the $X_i$ pull-back to Cartier divisors,
then so is $(\Til M \times_M \Til M')^{\sim}\to M$, where $(\Til M
\times_M \Til M')^{\sim}$ is the component of the fiber product
dominating $M$. So we may assume $\Til M'$ dominates~$\Til M$:
$$\xymatrix{
{\Til M'} \ar[r]_p \ar@/^1pc/[rr]^{\pi'} & {\Til M} \ar[r]_{\pi} & {M}
}\quad.$$
Denoting by a $'$ the corresponding subschemes in $\Til M'$, we have
then
$$\overline Y'=p^* \overline Y\quad,\quad \overline X'_i=p^* \overline
X_i\quad,\quad R'_i=p^* R_i\quad,$$
and the stated independence follows immediately.
\end{proof}

Before proving our inclusion-exclusion formula, we comment on
relations between the class we introduced and ordinary Segre classes.

\begin{prop}\label{elem}
\begin{enumerate}
\item\label{r1} For $r=1$: $s(Y;X;M)=s(X,M)$ is the Segre class of $X$
  in~$M$ (in particular, it is independent of $Y$).
\item\label{empt} If $Y=\emptyset$, then
  $s(Y;X_1,\dots,X_s;M)=s(X_1\cup\cdots\cup X_s,M)$, where
  $X_1\cup\cdots\cup X_s$ denotes the subscheme with ideal given by
  the product of the ideals of $X_1,\dots,X_s$.
\item $s(Y;X_1,\dots,X_s;M)$ is a birational invariant: if $f: M'\to
  M$   is a proper birational morphism, and $Y'=f^{-1}(Y)$,
  $X'_i=f^{-1}(X_i)$ (scheme-theoretically) then
$$f_* s(Y';X'_1,\dots,X'_s;M')=s(Y;X_1,\dots,X_s;M)\quad.$$
\end{enumerate}
\end{prop}

\begin{proof}
The first two statements follow immediately from the definition.
For the third: if $\pi':\Til M' \to M'$ is a resolving morphism for
$Y'$ and $X'_1$, \dots, $X'_s$, then $\pi'\circ f$ is resolving for
$Y$ and $X_1, \dots, X_s$, and the result follows easily.
\end{proof}

The second statement of Proposition~\ref{elem} is a property shared
with conventional Segre classes, cf.~Proposition~4.2 in
\cite{MR85k:14004}; it can be useful for concrete computations.

Here is the promised inclusion-exclusion principle:

\begin{theorem}[Inclusion--Exclusion formula]\label{main}
With notations as above, let\newline $Y\subset X_1\cap\dots\cap X_r$, and
assume $R_1\cdot\dots\cdot R_r=0$. Then
$$s(Y,M)=\sum_{s=1}^r (-1)^{s-1}\sum_{i_1<\dots<i_s} s(Y;
X_{i_1},\dots, X_{i_s};M)$$
in $A_*(X_1\cup\dots\cup X_r)$.
\end{theorem}

\begin{remark}\label{cases}
The condition on the $R_i$ seems hard to check, but it is
automatically verified in at least two situations:
\begin{enumerate}
\item if $Y$ in fact {\em equals\/} $X_1\cap\dots\cap X_r$
(scheme--theoretically); and
\item if $r>\dim M$.
\end{enumerate}
The first case is the one which best agrees with the set--theoretic
analogy; the second says (maybe surprisingly) that we could, if we
want, choose {\em all $X_i$ to be the same;\/} so long as they contain
$Y$, and we are choosing enough of them, then the formula must
hold. We will illustrate this point in the beginning of \S\ref{SM}.
\end{remark}

\begin{proof}
By Proposition~\ref{elem}, we can replace $M$ by a birational variety
obtained by a resolving morphism. Thus we may assume that $Y$ and the
$X_i$ are Cartier divisors, and we are reduced to proving that 
$$\frac{[Y]}{1+Y}-\sum_{s=1}^r (-1)^{s-1}\sum_{i_1<\dots<i_s}
\frac{[R_{i_1}+\dots+R_{i_s}+Y]}{1+R_{i_1}+\dots+R_{i_s}+Y}$$
is zero under the stated conditions; that is, that
$$\sum_{s=0}^r (-1)^s\sum_{i_1<\dots<i_s}
\frac{[R_{i_1}+\dots+R_{i_s}+Y]}{1+R_{i_1}+\dots+R_{i_s}+Y}$$
is 0 if $R_1\cdot\dots\cdot R_r=0$. But for this it is enough to show
that this last expression, as a formal power series in the $R_i$'s and
$Y$, is divisible by $R_1\cdot\dots\cdot R_r$; and in turn for this it
is enough to show that it is divisible by $R_r$; and finally, for this
it suffices to show that the expression vanishes when $R_r$ is set
to~0. Rewriting the expression by separating the part which contains
$R_r$ transforms it into
$$\sum_{s=0}^r (-1)^{s+1}\sum_{i_1<\dots<i_s<r}
\frac{[R_{i_1}+\dots+R_{i_s}+R_r+Y]}{1+R_{i_1}+\dots+R_{i_s}+R_r+Y}$$
plus
$$\sum_{s=0}^r (-1)^s\sum_{i_1<\dots<i_s<r}
\frac{[R_{i_1}+\dots+R_{i_s}+Y]}{1+R_{i_1}+\dots+R_{i_s}+Y}\quad.$$
Setting $R_r=0$ in the first summands makes it equal and opposite
in sign to the second; so the sum is 0 when $R_r=0$, as needed.
\end{proof}

Without hypotheses on the residual $R_i$ (via a resolving morphism
$\pi$) and if $r\le n=\dim M$, it is not hard to still say something
rather precise regarding the right-hand-side of our
inclusion-exclusion formula. As an illustration, consider the case
$r=n$; then the proof of the theorem shows that the difference between
the inclusion--exclusion formula and the Segre class of $Y$ is a
multiple of $\pi_*(R_1\cdot\dots\cdot R_n)$. This multiple can be
evaluated easily.

\begin{prop}\label{fact}
With notations as above, let $Y\subset X_1\cap\dots\cap X_n$, with
$n=\dim M$; then
$$\sum_{s=1}^r (-1)^{s-1}\sum_{i_1<\dots<i_s} s(Y;
X_{i_1},\dots, X_{i_s};M)=s(Y,M)+n!\, \pi_*(R_1\cdot \dots \cdot R_n)$$
in $A_*(X_1\cup\dots\cup X_r)$.
\end{prop}

\begin{proof}
Again replace $M$ by a resolving variety. The coefficient of $R_1\cdot
\dots\cdot R_n$ can be computed after assuming $Y=0$ and all $R_i$ are
equal to a fixed class $R$. 

With these positions the left-hand-side becomes 
$$\sum_{s=1}^n (-1)^{s-1} \binom ns \frac{s[R]}{1+sR}\quad;$$
the coefficient of $R^n$ in this expression is
$$\sum_{s=0}^n (-1)^{n-s} \binom ns s^n\quad.$$
So the result is implied by the combinatorial identity which follows.
\end{proof}

\begin{lemma}
$$\sum_{s=0}^n (-1)^{n-s} \binom ns s^n=n!$$
\end{lemma}

This is easily verified: $n^n$ counts the number of $n$--digit strings
from $n$ objects; $-n\cdot (n-1)^n$ takes away those which only use
$n-1$ objects; but those that used $n-2$ or less have been taken away
twice, so $+\binom n2 \cdot (n-2)^n$ evens out the count; except that
those that used $n-3$ or less have been counted too many times;
etc.~etc. The conclusion is that the formula in the lemma computes the
number of strings of length $n$ using exactly $n$ objects, that is,
$n!$.

As a side remark, we note that the argument in the proof of the
theorem shows, by the same token, that
$$\sum_{s=0}^n (-1)^{r-s} \binom ns s^r=0\quad\hbox{for $0\le r<n$}\quad;$$ 
which is also immediate from the counting argument just provided,
since there are {\em no\/} strings of length $n$ using fewer than $n$
objects.


\section{Applications and examples}\label{SM}

There are two main situations in which Theorem~\ref{main} applies
(cf.~Remark~\ref{cases}): when the schemes $X_i$ cut out $Y$,
and when there simply are enough of them.

To illustrate the second situation, let
$$s(Y;X^{(r)};M):=s(Y;\overbrace{X,\dots,X}^r;M)$$
where $X$ is any scheme containing $Y$. The classes
$$\sum_{s=1}^r (-1)^{s-1}\binom rs s(Y;X^{(s)};M)$$
may be seen as `successive approximations' of $s(Y,M)$: indeed, by
Proposition~\ref{elem}, part~\ref{r1}., and Theorem~\ref{main}, this
class equals $s(X,M)$ for $r=1$, and $s(Y,M)$ for $r\gg 0$.

\begin{example}
Let $M=\Pbb^5$, $X=\Pbb^4\subset M$, and $Y=$ a quadric surface in
a $\Pbb^3$ contained in $X$. The successive approximations of $s(Y,M)$
are (after push-forward to $M=\Pbb^5$):
\begin{alignat*}{5}
s(X,M) = &[\Pbb^4]&-&[\Pbb^3]&+&[\Pbb^2]&-&[\Pbb^1]&+&[\Pbb^0] \\
s(Y;X^{(2)};M) = && 2&[\Pbb^3]&-4&[\Pbb^2]&+6&[\Pbb^1]&-8&[\Pbb^0]\\
s(Y;X^{(3)};M) = && && -4&[\Pbb^2]& +4 &[\Pbb^1]& -8 &[\Pbb^0]\\
s(Y;X^{(4)};M) = && &&  2&[\Pbb^2]& +16&[\Pbb^1]& +22&[\Pbb^0]\\
s(Y;X^{(5)};M) = && &&  2&[\Pbb^2]& -8 &[\Pbb^1]& -98&[\Pbb^0]\\
s(Y;X^{(6)};M) = && &&  2&[\Pbb^2]& -8 &[\Pbb^1]& +22&[\Pbb^0]\\
\end{alignat*}
and they stabilize from this point on, having reached $s(Y,M)$.
\end{example}

Since successive approximations must eventually agree, we obtain a
relation recursively computing $s(Y;X^{(r)};M)$ for large
$r$. Explicitly:
$$s(Y;X^{(r+1)};M)=\sum_{s=1}^r (-1)^{r-s} \binom r{s-1}
s(Y;X^{(s)};M)\quad\text{for $r>\dim M$.}$$
More generally, Theorem~\ref{main} could be used to obtain a recursion
to compute the class $s(Y;X_1,\dots,X_r;M)$ for large $r$, in terms of
classes involving fewer of the $X_i$'s.

For $Y=\emptyset$ and all $X_i$ equal to a fixed hypersurface $X$, the
recursion says the following. Let $M$ be a nonsingular variety, and
let $\mathcal L$ be a line bundle on $M$. Denote by $c(rX)\in A_*M$
the total Chern class of a nonsingular section of $\mathcal L^{\otimes
  r}$ (if any exists). Then
\begin{itemize}
\item For $r>\dim M$:
$$c(rX)=\sum_{s=1}^{r-1} (-1)^{r-s-1} \binom rs c(s X)$$
\item For $r=\dim M$:
$$c(rX)=\sum_{s=1}^{r-1} (-1)^{r-s-1} \binom rs c(s X)-
(\dim M)!\, c_1(\mathcal L^\vee)^{\dim M}\cap [M]$$
\end{itemize}

The first formula is immediate from the above considerations and
Proposition~\ref{elem}, part~\ref{empt}; the second one follows
likewise by using Proposition~\ref{fact}. 

These facts are essentially obvious independently of the
results in this paper (exercise), but may deserve to be better known
then they seem to be. For instance, this recursive behavior applies to
the {\em Euler characteristic\/} of nonsingular sections of higher and
higher powers of a line-bundle. 

\begin{example}
The simplest possible example of this phenomenon has to be the case of
curves in $\Pbb^2$, where it can of course be immediately verified
from the genus formula. Even in this case, however, the recursive
recipe seems rather pretty: taking $\mathcal L=\mathcal O(1)$,
the successive Euler characteristics of degree $r$ curves, for
$r=1,2,3,\dots$ are
$$\bf 2,2,0,-4,-10,-18,-28,-40,\dots\quad;$$
in this case, the observation amounts to the fact that
$$\gathered
\tbinom 21\cdot \mathbf 2-2!\cdot 1=\mathbf 2\\
-\tbinom 31\cdot 2+\tbinom 32\cdot \mathbf 2= \mathbf 0\\
\tbinom 41\cdot 2-\tbinom 42\cdot 2+\tbinom 43\cdot \mathbf 0=\mathbf{-4} \\
-\tbinom 51\cdot 2+\tbinom 52\cdot 2-\tbinom 53\cdot 0+\tbinom 54
\cdot(\mathbf{-4})=\mathbf{-10}\\
\tbinom 61\cdot 2-\tbinom 62\cdot 2+\tbinom 63\cdot 0-\tbinom 64\cdot
(-4)+ \tbinom 65\cdot(\mathbf{-10})=\mathbf{-18}\\
\endgathered$$
etc.
\end{example} 

By the above considerations, the same structure governs all sequences
of Euler characteristics of successive powers of a line bundle over
every variety; if the dimension of $M$ is $n$, the structure is
visible from the $n$-th row of Pascal's triangle onward. The same
applies to the terms in all dimensions of the classes $c(rX)$.\vskip 6pt

A more substantial application of Theorem~\ref{main} is to the
computation of {\em SM-Segre classes,\/} that is (equivalently) of
Chern-Schwartz-MacPherson classes. The reader is addressed to
\cite{incexcI} and references therein for more thorough discussions of
these notions. As a quick reminder, the Chern-Schwartz-MacPherson
classes are classes defined in the Chow group of a possibly singular
variety $V$, agreeing with $c(TV)\cap [V]$ when $V$ is nonsingular,
and fitting a functorial prescription (originally envisioned by
Deligne and Grothendieck; see \cite{MR91h:14010} for a nice
treatment). There has been substantial activity in recent years,
comparing these classes to other classes generalizing the classes of
the tangent bundle. These other classes are typically defined in terms
of a Segre class; thus it is natural to ask whether the
Chern-Schwartz-MacPherson classes admit a description in terms of
Segre classes.

This is indeed the case for hypersurfaces (\cite{MR2001i:14009}): the
difference between Chern-Schwartz-MacPherson classes and the classes
of the virtual tangent bundle of a hypersurface $X$ can be measured
effectively in terms of the Segre classes of the {\em singularity
  subscheme\/} of $X$. This difference has been named {\em Milnor
  class,\/} and studied in a different context; see \cite{MR1795550}
for the hypersurface case, and \cite{math.AG/0202175} for the state of
the art on this point of view (now embracing a substantially more
general class of varieties).

This otherwise substantial progress does not seem to yield a general
expression for the Chern-Schwartz-MacPherson class in terms of Segre
classes. Theorem~\ref{relSM} below is a (small) step in this direction. 

Let $M$ be a {\em nonsingular\/} ambient variety.

\begin{defin}
We say that a reduced subscheme $X\subset M$ is {\em almost
  nonsingular\/} if its irreducible components $X_i$ are nonsingular,
  and $X_i\cap X_j=Y$ (as schemes) for all $i\ne j$, where $Y$ is a
  nonsingular variety independent of $i$, $j$.
\end{defin}

\begin{example}
The simplest example of an almost nonsingular variety is the
transversal union of two nonsingular hypersurfaces in $M$. Another
typical example would be the union of any number of lines through a
point in a projective space, or more generally the union of a
collection of general subspaces of dimension $\le \frac n2$ containing
a fixed subspace $Y$ in $\Pbb^n$.

In particular, note that {\em almost nonsingular\/} varieties are {\em
  not} necessarily pure-dimensional, Cohen-Macaulay, complete
intersections\dots
\end{example}

\begin{theorem}\label{relSM}
Let $X\subset M$ be almost nonsingular, let $X_1,\dots,X_r$ be its
components, and let $Y$ be the common pairwise intersection of the
components of $X$. Then
$$\csm(X)=c(TM)\cap s(Y;X_1,\dots,X_r;M)\quad.$$
\end{theorem}

\begin{proof}
Induction on the number $r$ of components of $X$. If $r=1$, then $X$
is nonsingular, so $\csm(X)=c(TX)\cap [X]=c(TM)\cap s(X,M)$; so the
formula holds, by Proposition~\ref{elem}, part~\ref{r1}.

For $r>1$, observe that since $\csm$ itself satifies
inclusion-exclusion (as a particular instance of functoriality) then
$$\csm(Y)=\sum_{s=1}^{r-1} (-1)^{s-1}\sum_{i_1<\dots<i_s}
\csm(X_{i_1}\cup\dots\cup X_{i_s})+(-1)^{r-1}\csm(X)\quad.$$
Capping by $c(TM)^{-1}$, and using the induction hypothesis:
$$s(Y,M)=\sum_{s=1}^{r-1} (-1)^{s-1}\sum_{i_1<\dots<i_s}
s(Y;X_{i_1},\dots,X_{i_s};M) + (-1)^{r-1}c(TM)^{-1}\cap \csm(X)$$
The result follows by comparing this formula with Theorem~\ref{main}.
\end{proof}

Using the definition of {\em SM-Segre class\/} $\ssm(X,M)$ introduced
in \cite{incexcI}:

\begin{corol}\label{surpr}
With notations as in Theorem~\ref{relSM},
$$\ssm(X,M)=s(Y;X_1,\dots,X_r;M)\quad.$$
\end{corol}

\begin{proof}
This now follows from Theorem~\ref{relSM} and \cite{incexcI},
Theorem~3.1.
\end{proof}

Corollary~\ref{surpr} gives a suprisingly straightforward computation
of $\ssm(X,M)$ in the particular case of almost nonsingular
varieties. The definition of $\ssm(X,M)$ given in \cite{incexcI} is a
certain combination of the Segre classes of hypersurfaces cutting out
$X$ in $M$ and of their singularity subschemes; it would seem very
difficult to have any control over it in that form. The expression in
terms of a union Segre class is by contrast very direct, and matches
well the intuition (from the case of a hypersurface,
cf.~\cite{MR96d:14004}, \S1) that the right class should be
obtained by `removing' a copy of the singularity subscheme.

\begin{example}\label{repr}
This is perhaps the main punch-line in the whole paper, so we expand
on this point. As mentioned above, the simplest example of an almost
nonsingular variety is the union $X$ of two transversal, nonsingular
hypersurfaces $X_1$, $X_2$. If $(F_1)$, $(F_2)$ are (local) ideals
for $X_1$, $X_2$, then $X$ has ideal $(F_1F_2)$, and hence
singularity subscheme $Y$ with ideal $(F_1F_2,F_1\,
dF_2+F_2\,dF_1)$. As $X_1$ and $X_2$ are transversal, $dF_1$ and
$dF_2$ are linearly independent at every point of $Y$; hence the ideal
of $Y$ is simply $(F_1,F_2)$: that is, $Y=X_1\cap X_2$.

Now using Theorem~\ref{relSM} we get
$$\csm(X)=c(TM)\cap s(Y;X_1,X_2;M)\quad.$$
Blowing up along $Y$, we can write $s(Y;X_1,X_2;M)$ as the
push-forward of 
$$\frac{[X_1+X_2-E]}{1+X_1+X_2-E}\quad,$$
where $E$ denotes the exceptional divisor and we use $X$ for the
inverse image of its namesake in $M$; in this sense we `remove one
copy of $Y$ from $X$'. Simple manipulations (using the notational
device from \cite{MR2001i:14009}, \S1.4) give
$$\frac{[X-E]}{1+X-E} =\frac {[X]}{1+X}+\frac 1{1+X}\cdot \frac
     {-[E]}{1+X-E}
= \frac{[X]}{1+X}+\frac 1{1+X}
     \left(\left(\frac{[E]}{1+E}\right)^\vee 
\otimes \mathcal O(X)\right),$$
and finally pushing forward we get
$$s(Y;X_1,X_2;M)=s(X,M)+c(\mathcal O(X))^{-1}\cap
\left(s(Y,M)^\vee\otimes \mathcal O(X)\right)\quad.$$
This is precisely the formula giving $s(X\setminus Y,M)$ (cf.~Lemma
I.3 in \cite{MR2001i:14009}). That is, although the case of `almost
nonsingular varieties' is extremely special, it is representative
enough to capture precisely the correct definition for the class
studied in full generality in \cite{MR2001i:14009}. 
\end{example}

We hope that studying union Segre classes for similarly representative
cases in higher codimension may suggest extensions of the result of
\cite{MR2001i:14009}. 

Conversely, the remarkable properties of $\ssm(X,M)$ (see
\cite{incexcI}, Theorem~2.3) tell us something about the union Segre
class, at least in the case of almost nonsingular varieties. For
instance: 

\begin{corol}\label{indep}
If $X$ is almost nonsingular, with components $X_i$, and $Y=X_i\cap
X_j$ for $i\ne j$, and $X\subset M$, $X\subset M'$ are embeddings of
$X$ in nonsingular varieties, then
$$s(Y;X_1,\dots,X_r;M)=c(TM)^{-1}c(TM')\cap
s(Y;X_1,\dots,X_r;M') \quad.$$
\end{corol}

\begin{proof}
By Theorem~\ref{relSM}, both classes $c(TM)\cap
s(Y;X_1,\dots,X_r;M)$ and $c(TM')\cap s(Y;X_1,\dots,X_r;M')$
compute $\csm(X_1\cup\cdots \cup X_r)$.
\end{proof}

Such observations indicate that union Segre classes are better behaved
than one may expect. For example, they must be less sensitive to
scheme structure than ordinary Segre classes.

\begin{example}
The {\em ordinary\/} Segre class of the union $X$ of three distinct
lines $L_1$, $L_2$, $L_3$ through a point $p$ in $\Pbb^3$ depends on
the position of the lines.
$$\epsffile{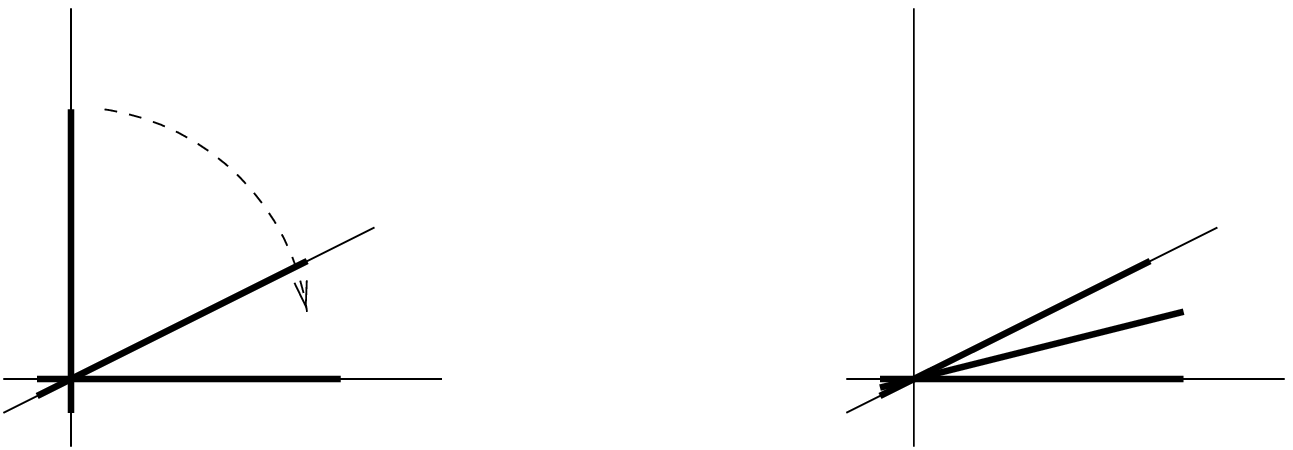}$$
If the three lines are {\em not\/} coplanar, then a direct computation
shows that $s(X,\Pbb^3)=[X]-10[\Pbb^0]$; if the lines {\em are\/}
coplanar, then $s(X,\Pbb^3)=[X]-12[\Pbb^0]$. Here we are taking the
reduced structure on the union, and in a sense this is the problem: as
one of the lines moves to become coplanar with the others, the flat
limit acquires an embedded point at the origin; in other words, the
reduced union of coplanar lines is not a limit of the union of
configurations of noncoplanar ones.

The union Segre class does not detect this difference. Using
Definition~\ref{usc}, or Theorem~\ref{relSM}, one computes that 
$$s(p;L_1,L_2,L_3;\Pbb^3)=[X]-12[\Pbb^0]$$
regardless of the directions of the lines, so long as these remain
distinct.
\end{example}

On the other hand, properties such as Corollary~\ref{indep} do depend
on rather subtle information, in that they do not extend blindly to
all union Segre classes.

\begin{example}
Consider the union $X$ of two distinct lines $L_1$, $L_2$ through a
point, with ideal given by the {\em product\/} of the ideals of the
components. 
$$\epsffile{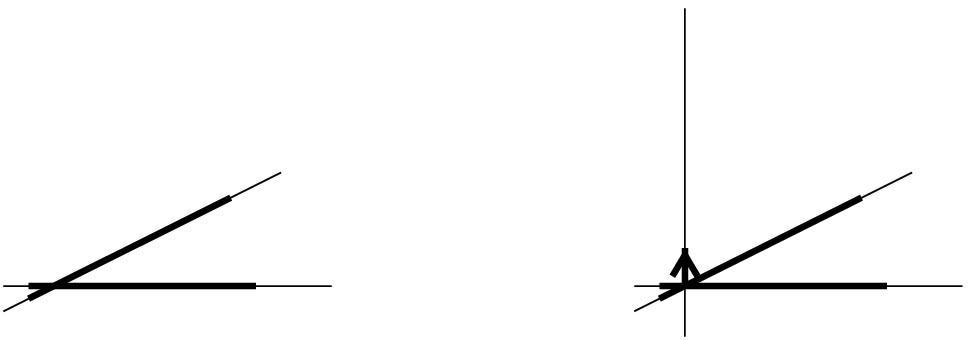}$$
This scheme depends on the ambient variety: in $\Pbb^2$ it is reduced,
in $\Pbb^3$ it acquires an embedded component supported at the
intersection point. Computing the ordinary Segre classes, and applying
Proposition~\ref{elem}, part~\ref{empt}., gives
$$s(\emptyset;L_1,L_2;\Pbb^2)=[X]-4[\Pbb^0]=
s(\emptyset;L_1,L_2;\Pbb^3)\quad;$$
thus
$$s(\emptyset;L_1,L_2;\Pbb^2)\ne c(T\Pbb^2)^{-1}c(T\Pbb^3)\cap
s(\emptyset;L_1,L_2;\Pbb^3)=[X]-6[\Pbb^0]\quad.$$
However, by Corollary~\ref{indep} we must have
$$s(p;L_1,L_2;\Pbb^2)= c(T\Pbb^2)^{-1}c(T\Pbb^3)\cap
s(p;L_1,L_2;\Pbb^3)\quad,$$
where $p$ is the intersection point. Using Definition~\ref{usc}, the
reader can check `by hand' that $s(p;L_1,L_2;\Pbb^2)=[X]-3[\Pbb^0]$,
$s(p;L_1,L_2;\Pbb^3)=[X]-5[\Pbb^0]$, in agreement with this formula.
\end{example}

Summarizing, the inclusion-exclusion principle draws a connection
between the union Segre classes considered here and the SM-Segre
classes considered in \cite{incexcI}. In view of applications to the
study of Chern-Schwartz-MacPherson and Milnor classes, it would be
very interesting to establish the precise circumstances under which
union and SM-Segre classes agree.


\end{document}